\documentclass{lms}

\usepackage{graphicx,color,epsf}
\usepackage [cmtip,arrow]{xy}
\xyoption{all}
\usepackage {pb-diagram,pb-xy} 
\usepackage{amsmath,amsfonts}

\newtheorem{theorem}{Theorem}[section] 
\newtheorem{lemma}[theorem]{Lemma}     
\newtheorem{corollary}[theorem]{Corollary}
\newtheorem{proposition}[theorem]{Proposition}

\newnumbered{assertion}{Assertion}    
\newnumbered{conjecture}{Conjecture}  
\newnumbered{definition}{Definition}
\newnumbered{hypothesis}{Hypothesis}
\newnumbered{remark}{Remark}
\newnumbered{note}{Note}
\newnumbered{observation}{Observation}
\newnumbered{problem}{Problem}
\newnumbered{question}{Question}
\newnumbered{algorithm}{Algorithm}
\newnumbered{example}{Example}
\newunnumbered{notation}{Notation} 

\newcommand{\Real}{{\mathbb R}}

\newcommand{\R}{{\rm  R}}

\newcommand{\eps}{\varepsilon}

\newcommand{\x}{\mathbf{x}}

\newcommand{\y}{\mathbf{y}}

\newcommand{\z}{\mathbf{z}}

\newcommand{\A}
 {\mathcal{A}}

\newcommand {\RR} {{\mathcal R}}

\newcommand {\hide}[1]{}

\newcommand {\CT} {{\rm CT}}
\newcommand {\OT} {{\rm OT}}
\newcommand {\BT} {{\rm BT}}

\newcommand {\HH} {{\rm H}}

\newcommand {\Ann} {{\rm Ann}}

\renewcommand{\a}{\mathbf{a}}

\renewcommand{\x}{\mathbf{x}}

\renewcommand{\y}{\mathbf{y}}

\renewcommand{\z}{\mathbf{z}}

\title[Combinatorial complexity in o-minimal geometry]%
{Combinatorial complexity in o-minimal geometry
}

\author{Saugata Basu}

\classno{52C45}

\extraline{
The author was supported in part by an NSF 
grant CCF-0634907. 
}
\begin{document}
\maketitle

\begin{abstract}
In this paper we prove tight bounds on the combinatorial and topological 
complexity
of sets defined in terms of $n$ definable sets belonging to some fixed
definable family of sets in an o-minimal structure. This generalizes the
combinatorial parts of similar bounds known in the case of semi-algebraic and 
semi-Pfaffian sets, and as a result vastly increases the applicability
of results on combinatorial and topological complexity of
arrangements studied in discrete and computational geometry. As a sample
application, 
we extend a Ramsey-type theorem due to Alon et al. \cite{APPRS}, originally 
proved for semi-algebraic sets of fixed description complexity 
to this more general setting.
\end{abstract}
\maketitle

\section{Introduction}
\label{sec:intro}
Over the last twenty years there has been a lot of work on bounding the
topological complexity (measured in terms of their Betti numbers)
of several different classes of subsets of $\Real^k$ -- most notably
semi-algebraic and semi-Pfaffian sets. The usual setting for
proving these bounds is  as follows. One considers a semi-algebraic 
(or semi-Pfaffian)
set $S \subset \Real^k$ defined by a Boolean formula
whose atoms consists of $P >0, P = 0, P < 0, \; P \in {\mathcal P},$
where ${\mathcal P}$ is a set of polynomials (resp. Pfaffian 
functions) of degrees bounded by a parameter (resp. whose
Pfaffian complexity is bounded by certain parameters) and 
$\#{\mathcal P} = n.$ 
It is possible to obtain  bounds on the Betti numbers of
$S$ in terms of $n,k$ and the parameters bounding the complexity of the
functions in ${\mathcal P}$.

\subsection{Known bounds in the semi-algebraic and semi-Pfaffian cases}
In the semi-algebraic case, if we assume that the degrees of the polynomials
in ${\mathcal P}$ are bounded by $d$, 
and denoting by $b_i(S)$ the $i$-th Betti number of $S$,
then it is shown in \cite{GaV} that,

\begin{equation}
\label{eqn:bound1}
\sum_{i \geq 0} b_i(S)  \leq  n^{2k} O(d)^k. 
\end{equation}
 
A similar bound is also shown for semi-Pfaffian sets \cite{GaV}.

In another direction, we also have reasonably tight bounds on the sum of 
the Betti numbers
of the realizations of all realizable sign conditions of the family
${\mathcal P}$. A sign condition on ${\mathcal P}$ is an element of
$\{0,1,-1\}^{{\mathcal P}}$, and the realization of a sign condition 
$\sigma$ is the set,
\[
\RR(\sigma) = \{ \x \in \R^k \; \mid \; {\rm sign}(P(\x)) = \sigma(P), \forall
P \in {\mathcal P}\}.
\]
It is shown in \cite{BPR02} that,

\begin{equation}
\label{eqn:bound2}
\sum_{\sigma \in \{0,1,-1\}^{{\mathcal P}}} 
        b_i(\RR(\sigma)) \leq \sum_{j=0}^{k - i} 
              {n \choose j} 4^{j}  d(2d-1)^{k-1} 
 =   n^{k-i} O(d)^k.
\end{equation}

We refer the reader to \cite{B00,BPR02,GVsurvey,GaV,GVZ04}, as well as 
the survey article \cite{BPR10} for a comprehensive history of the
work leading up to the above results, as well as several other interesting
results in this area.

\subsection{Combinatorial and algebraic complexity}
Notice that the bounds in (\ref{eqn:bound1}) and 
(\ref{eqn:bound2}) are products of two quantities -- 
one that depends only on $n$ (and $k$),
and another part which is independent of $n$,  but
depends on the parameters controlling the
complexity of individual elements of ${\mathcal P}$ (such as degrees
of polynomials in the semi-algebraic case, 
or the  degrees and the length of the Pfaffian chain defining
the functions in the Pfaffian case). It is customary to
refer to the first part as the {\em combinatorial part}
of the complexity, and the latter as the algebraic (or Pfaffian)  part. 
Moreover, the algebraic or the Pfaffian parts of the bound depend on 
results whose proofs  involve  Morse theory (for instance, the well known 
Oleinik-Petrovsky-Thom-Milnor bounds on the Betti numbers of 
real varieties \cite{OP,T,Milnor2}).

While understanding the algebraic part of the complexity 
is a very important problem, in several
applications, most notably in discrete and computational geometry, 
it is the combinatorial part of the complexity that is of primary interest
(the algebraic part is assumed to be bounded by a constant). 
The motivation behind this point of view is the following.
In problems in discrete and computational geometry,
one typically encounters arrangements of a large number of  
objects in $\Real^k$
(for some fixed $k$), where each object is of ``bounded description 
complexity'' 
(for example, defined by a polynomial inequality of degree
bounded by a constant). Thus, it is the number of objects that
constitutes the important parameter, and the algebraic complexity
of the individual objects are thought of as small constants.
It is this second setting that is our primary interest in this paper.

The main results of this paper generalize  (combinatorial parts of) 
the bounds in (\ref{eqn:bound1}) and (\ref{eqn:bound2}) to  
sets which are definable in an arbitrary 
o-minimal structure over a  real closed field $\R$ (see Section
\ref{sec:ominimal} below for the definition of an o-minimal structure
and definable sets).

Instead of only considering sets having 
``bounded description complexity'', we allow the sets in an arrangement
${\mathcal A}$  to
be fibers of some fixed definable map 
$\pi: T \rightarrow \R^{\ell}$, where
$T  \subset \R^{k + \ell}$ is a definable set.
This vastly expands the applicability of results 
concerning complexity of arrangements in discrete and computational
geometry, 
since it is no longer necessary that the objects in the arrangements
be defined only in terms of polynomials. 
As we will see shortly, the sets we consider are allowed to be
fairly arbitrary.
They include sets defined by restricted analytic functions, 
including (but not by any means restricted to) polynomials,
Pfaffian functions such as exponential, logarithmic, trigonometric, inverse
trigonometric functions, subject to some mild conditions. 
All hitherto considered families of objects in the computational geometry 
literature, such as hyperplanes, simplices, and more generally
sets having  bounded description complexity are  special
instances of this general definition.
We also consider sets
belonging to the Boolean algebra generated by 
$n$ sets in $\R^k$ each of which is a fiber of a fixed definable map.
We prove tight bounds on the Betti numbers, 
the topological complexity of projections,
as well as on the complexity of cylindrical decomposition of such sets, 
in terms of $n$ and $k$. The role of the algebraic complexity is played by
a constant that depends only on the particular definable family. In this
way, we are able to generalize the notion of combinatorial complexity
to definable sets over an arbitrary o-minimal structure. 

Apart from the intrinsic mathematical interest of the results proved in the
paper, we believe that the techniques used to prove them would be of 
interest to researchers in discrete and computational geometry. 
We show that most (if not all) results
on the complexity of arrangements are consequences of a set of very simple
and well-studied axioms (those defining o-minimal structures). 
Many widely used techniques in the study of arrangements are strongly 
dependent on the assumption that the sets under consideration 
are semi-algebraic. For example,
it is common to consider  real algebraic varieties of fixed degree as 
hyperplane sections of the corresponding Veronese variety in a higher
(but still fixed) dimensional space -- a technique called ``linearization''
in computational geometry literature (see \cite{Agarwal}). 
Obviously, such methods fail if the given sets are not semi-algebraic. 
Our methods make no use of semi-algebraicity of the objects, 
nor bounds derived from Morse theory
such as the classical Oleinik-Petrovsky-Thom-Milnor bounds on 
Betti numbers of real algebraic varieties. 
We believe that this point of view simplifies proofs, 
and simultaneously generalizes vastly the class of objects which are allowed, 
at the same time 
getting rid of unnecessary assumptions such as requiring the 
objects to be in  general position. It is likely  that the techniques 
developed here will find further applications
in the combinatorial study of arrangements other than those 
discussed in this paper.

\subsection{Arrangements in computational geometry}
We now make precise the notions of arrangements, cells and their complexities,
following their usual definitions in discrete and computational
geometry \cite{Agarwal,Matousek}.

Let
${\mathcal A} = \{S_1,\ldots,S_n\}$, such that each $S_i$ is a
subset of $\R^k$ belonging to some ``simple'' class of sets. 
(We will define the class of admissible sets that we consider 
precisely in Section \ref{subsec:admissible} below). 

For $I \subset \{1,\ldots,n\}$, we let ${\mathcal A}(I)$ denote the set
\[
\bigcap_{i \in I \subset [1\ldots n]} S_i  \;\; \cap 
\bigcap_{j \in [1\ldots n]\setminus I}
\R^k \setminus S_j, 
\]
and it is customary to call a connected component of $A(I)$ a
\emph{cell} of the arrangement (even though it might not be a cell in the sense
of topology).  We let ${\mathcal C}({\mathcal A})$ denote the set 
of all non-empty cells of the arrangement ${\mathcal A}$.

The cardinality of ${\mathcal C}({\mathcal A})$ is called the 
{\em combinatorial complexity} of the arrangement ${\mathcal A}$. 
Since different 
cells of an arrangement might differ topologically, 
it makes sense to give more weight to a topologically complicated
cell than to a topologically simple one in the definition of complexity.
With this in mind we define  (following \cite{Basu3}) the {\em topological
complexity} of a cell to be the sum of its Betti numbers (the ranks of 
singular homology groups of the cell).

The class of sets  usually considered in the study of arrangements
are sets with ``bounded description complexity'' (see \cite{Agarwal}). 
This means that each set in the arrangement is defined by a first order 
formula in the language of ordered fields involving 
at most a constant number polynomials whose  degrees are also 
bounded by a constant.
Additionally, there is often a requirement that the sets be in ``general
position''. The precise definition of ``general position'' varies with
context, but often involves restrictions such as:  the sets in the 
arrangements  are smooth manifolds, intersecting transversally.

\subsection{Arrangements over an o-minimal structure}
O-minimal structures present a 
natural mathematical framework to state and prove results on
the complexity of arrangements. 
In this paper  we consider arrangements whose members come from some fixed
definable family in an o-minimal structure (see below for definitions).
The usual notion of ``bounded description complexity'' turns out to be
a special case of this more general definition.

\subsubsection{O-minimal structures}
\label{sec:ominimal}
O-minimal structures were invented and first studied by
Pillay and Steinhorn in the pioneering papers
\cite{PS1,PS2} 
in part to show that the tame topological
properties exhibited by the class of semi-algebraic sets are consequences
of a set of few simple axioms. 
Later the theory was further
developed through contributions of other researchers, most notably
van den Dries, Wilkie, Rolin, Speissegger amongst others
\cite{Dries2,Dries3,Dries4,Wilkie,Wilkie2,Rolin}. We particularly
recommend the book by van den Dries \cite{Dries} and the notes by
Coste \cite{Michel2} for an easy introduction to the topic as well as the
proofs of the basic results that we use in this paper.

An o-minimal structure on a real closed field $\R$ is just
a class of subsets of $\R^k, k\geq 0$, (called the 
{\em definable sets} in the structure) satisfying these axioms (see below).
The class of semi-algebraic sets is one obvious example of such a structure,
but in fact there are much richer classes of sets which have been proved
to be o-minimal (see below). For instance, subsets of $\Real^k$ 
defined in terms inequalities involving not just polynomials, but also
trigonometric and exponential functions on restricted domains  have been 
proved to be o-minimal.

We now formally define o-minimal structures (following \cite{Michel2}).    
\begin{definition}
\label{def:o-minimal}
An o-minimal structure on a real closed field $\R$ is a sequence 
${\mathcal S}(\R) = ({\mathcal S}_n)_{n \in {\mathbb N}}$, where 
each ${\mathcal S}_n$ is a collection of subsets of $\R^n$, satisfying the 
following axioms \cite{Michel2}. 

\begin{enumerate}
\item
All algebraic subsets of $\R^n$ are in ${\mathcal S}_n$.
\item
The class ${\mathcal S}_n$ is closed under complementation and
finite unions and intersections.
\item
If $A \in {\mathcal S}_m$ and $B \in {\mathcal S}_n$ then
$A \times B \in {\mathcal S}_{m+n}$.
\item
If $\pi: \R^{n+1} \rightarrow \R^{n}$ is the projection map on the
first $n$ co-ordinates and $A \in {\mathcal S}_{n+1}$, then 
$\pi(A)  \in {\mathcal S}_n$.
\item
The elements of ${\mathcal S}_1$ are precisely finite unions of points
and intervals.
\end{enumerate}
\end{definition}

\subsubsection{Examples of o-minimal structures}
\begin{example}
\label{example:sa}
Our first example of an o-minimal structure
${\mathcal S}(\R)$, is the o-minimal structure over a real closed field $\R$
where
each ${\mathcal S}_n$  is 
the class of semi-algebraic subsets of $\R^n$. 
It follows easily from the Tarski-Seidenberg principle (see \cite{BCR})
that the class of sets ${\mathcal S}_n$
satisfies the axioms in Definition \ref{def:o-minimal}. We will denote
this o-minimal structure by ${\mathcal S}_{{\rm sa}}(\R)$.
\end{example}

If Example \ref{example:sa} was the only example of o-minimal structure
available then the notion of o-minimality would not be very interesting. 
However, there are many more examples (see for example
\cite{Dries,Dries2,Dries3,Dries4,Rolin,Wilkie,Wilkie2}).

\begin{example} \cite{Wilkie}
Let ${\mathcal S}_n$ be the images in $\Real^n$ under the projection maps
$\Real^{n+k} \rightarrow \Real^n$ of sets of the form 
$\{(\x,\y) \in \Real^{n+k} \mid P(\x,\y,e^{\x},e^{\y}) = 0\}$, where $P$ is
a real polynomial in $2(n+k)$ variables, and 
$e^{\x} = (e^{x_1},\ldots,e^{x_n})$ and
$e^{\y} = (e^{y_1},\ldots,e^{y_k})$.
We will denote this o-minimal structure over $\Real$ 
by ${\mathcal S}_{{\rm exp}}(\Real)$.
\end{example}

\begin{example} \cite{Gabrielov96}
Let ${\mathcal S}_n$ be the images in $\Real^n$ under the projection maps
$\Real^{n+k} \rightarrow \Real^n$ of sets of the form 
$\{(\x,\y) \in \Real^{n+k} \mid P(\x,\y)=0\}$, where $P$ is
a restricted analytic function in $n+k$ variables.
A restricted analytic function in $N$ variables is an analytic function
defined on an open neighborhood of $[0,1]^N$ restricted to
$[0,1]^N$ (and extended by $0$ outside).
We will denote this o-minimal structure over $\Real$  
by ${\mathcal S}_{{\rm ana}}(\Real)$.
\end{example}
The o-minimality of the last two classes
are highly non-trivial theorems.

\subsection{Admissible sets}
\label{subsec:admissible}
We now define the sets that will play the role of objects of ``constant
description complexity'' in the rest of the paper.

\begin{definition}
\label{def:admissible}
Let ${\mathcal S}(\R)$ be an o-minimal structure on a real
closed field $\R$ and let $T \subset \R^{k+\ell}$
be a definable set. 
Let $\pi_1: \R^{k+\ell} \rightarrow \R^{k}$
(resp. $\pi_2: \R^{k+\ell}  \rightarrow \R^{\ell}$), 
be the projections onto the first $k$ (resp. last $\ell$) co-ordinates.

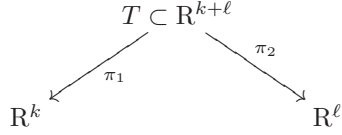
\[
\xymatrix{
& T \subset \R^{k+\ell} \ar[ld]^{\pi_1} \ar[rd]^{\pi_2} & \\
\R^{k} & & \R^{\ell}
}
\]

We will call a subset $S$ of $\R^k$ to be a $(T,\pi_1,\pi_2)$-set if
\[
S = \pi_1(\pi_2^{-1}(\y)\cap T)
\]
for some $\y \in \R^{\ell}$,
and when the context is clear we will denote 
$T_{\y} = \pi_1(\pi_2^{-1}(\y)\cap T)$.
In this paper, we will consider finite families of $(T,\pi_1,\pi_2)$-sets,
where $T$ is some fixed 
definable set for each such family,
and we will call a family of $(T,\pi_1,\pi_2)$-sets to be a 
$(T,\pi_1,\pi_2)$-family. We will also sometimes refer to a finite
$(T,\pi_1,\pi_2)$-family as an {\em arrangement} of $(T,\pi_1,\pi_2)$-sets.
\end{definition}

For any  definable set $X \subset \R^k$, 
we let $b_i(X)$ denote the $i$-th Betti number of $X$,
and we let $b(X)$ denote $\sum_{i\geq 0} b_i(X)$.
We define the {\em topological complexity}
of an  arrangement ${\mathcal A}$ of $(T,\pi_1,\pi_2)$-sets
to be the number
\[
\sum_{D \in {\mathcal C}({\mathcal A})} \sum_{i=0}^k b_i(D).
\]   

\begin{remark}
We remark here that
for o-minimal structures over an arbitrary real closed field $\R$, ordinary 
singular homology is not well defined. 
Even though o-minimal versions of singular
co-homology theory, as well \v{C}ech co-homology theory, 
has been developed recently (see \cite{EP,EW}),
in this paper we take a simpler approach
and use a modified homology theory 
(which agrees with singular homology in case $\R = {\mathbb R}$
and which is homotopy invariant)  as done in \cite{BPRbook2} in case
of semi-algebraic sets over arbitrary real closed fields (see 
\cite{BPRbook2}, page 279). The underlying idea behind that definition
is as follows. Since closed and bounded semi-algebraic (as well as
definable) sets are finitely triangulable, simplicial homology
is well defined for such sets. 
Furthermore, it is shown in \cite{BPRbettione} that it is possible to replace
an arbitrary semi-algebraic set by a closed and bounded one which is
homotopy equivalent to the original set.
We prove an analogous result for arbitrary definable sets in this paper
(see Theorem \ref{the:GV} below).  
We now define the 
homology groups of the original set to be the simplicial homology groups of
the closed and bounded definable set which is homotopy equivalent to it. It
is clear that this definition is homotopy invariant.
\end{remark}

We now give a few examples to show that arrangements of objects of
bounded description complexities are included in the class of arrangements
we study, but our class is much larger since $T$ need not be semi-algebraic.

\subsubsection{Examples}
\begin{example}
Let ${\mathcal S}(\R)$ be the o-minimal structure 
${\mathcal S}_{{\rm sa}}(\R)$. 
Let $T \subset \R^{2k+1}$ be the semi-algebraic
set defined by 
\[
T = \{(x_1,\ldots,x_k,a_1,\ldots,a_k,b) \mid \langle \a,\x\rangle -b = 0\}
\]
(where we denote $\a = (a_1,\ldots,a_k)$ and $\x = (x_1,\ldots,x_k)$),
and $\pi_1$ and $\pi_2$ are the projections onto the first $k$ and last $k+1$
co-ordinates respectively. A $(T,\pi_1,\pi_2)$-set is clearly a hyperplane
in $\R^k$ and vice versa.
\end{example}

\begin{example}
Again, let ${\mathcal S}(\R)$ be the o-minimal structure 
${\mathcal S}_{{\rm sa}}(\R)$. 
Let $T \subset \R^{k + k(k+1)}$ be the semi-algebraic
set defined by 
\[
T = \{(\x,\y_0,\ldots,\y_{k}) \mid \x,\y_0,\ldots,\y_{k} \in \R^k,
\x \in {\rm conv}(\y_0,\ldots,\y_{k})\},
\]
and $\pi_1$ and $\pi_2$ the projections onto the first $k$ and last $k(k+1)$
co-ordinates respectively. 
A $(T,\pi_1,\pi_2)$-set is a (possibly degenerate)
$k$-simplex in $\R^k$ and vice versa.
\end{example}

Arrangements of hyperplanes as well as simplices have been well studied in 
computational geometry, and thus the two previous examples do not
introduce anything new. We now discuss an example which could not be
handled by the existing techniques in computational geometry, such as
linearization. 

\begin{example}
Now, let ${\mathcal S}(\Real)$ be the o-minimal structure 
${\mathcal S}_{{\rm exp}}(\Real)$. 
Let $T \subset \Real^{k + m(k+1)}$ be the 
set defined by 
$$
\displaylines{
T = \{(\x,\y_1,\ldots,\y_{m},a_1,\ldots,a_m) \mid \x,\y_1,\ldots,\y_m \in 
\Real^k, a_1,\ldots,a_m \in \Real, \cr
x_1,\ldots,x_k >  0,  \sum_{i=0}^{m} a_i \x^{\y_i} = 0\},
}
$$
and $\pi_1: \Real^{k + m(k+1)} \rightarrow \Real^k$ and 
$\pi_2: \Real^{k + m(k+1)} \rightarrow \Real^{m(k+1)}$ be
the projections onto the first $k$ and last $m(k+1)$
co-ordinates respectively. 
It can be shown that $T$ is definable in the structure 
${\mathcal S}_{{\rm exp}}(\Real)$. 
The $(T,\pi_1,\pi_2)$-sets in this example include (amongst others)
all semi-algebraic sets
consisting of intersections with the positive  orthant
of all real algebraic sets defined by a polynomial
having at most $m$ monomials (different sets of monomials are allowed
to occur in different polynomials).
\end{example}

\begin{definition}
\label{def:A-sets}
Let ${\mathcal A} = \{S_1,\ldots,S_n\}$, such that each $S_i \subset \R^k$ 
is a $(T,\pi_1,\pi_2)$-set.
For $I \subset \{1,\ldots,n\}$, we let ${\mathcal A}(I)$ denote the set
\begin{equation}
\label{eqn:basic}
\bigcap_{i \in I \subset [1\ldots n]} S_i  \;\; \cap 
\bigcap_{j \in [1\ldots n]\setminus I}
\R^k \setminus S_j, 
\end{equation}
and we will call such a set to be a basic ${\mathcal A}$-set.
We will denote by ${\mathcal C}({\mathcal A})$
the set of non-empty connected components of all basic
${\mathcal A}$-sets.
 
We will call definable subsets $S \subset \R^k$ defined by a
Boolean formula whose atoms are of the form,
$x \in S_i, 1 \leq i \leq n$,  a ${\mathcal A}$-set. 
A ${\mathcal A}$-set is thus a union of basic ${\mathcal A}$-sets.
If  $T$ is closed, and the Boolean formula
defining $S$ has no negations, then $S$ is closed by definition
(since each $S_i$ being homeomorphic to the intersection of $T$
with a closed set $\pi^{-1}(\y)$ for some $\y \in\R^{\ell}$
is closed) and we call such a set an ${\mathcal A}$-closed set.

Moreover, if $V$ is any closed 
definable subset of $\R^k$,
and $S$ is an ${\mathcal A}$-set (resp. ${\mathcal A}$-closed set), then we
will call $S \cap V$  an $({\mathcal A},V)$-set 
(resp. $({\mathcal A},V)$-closed set).
\end{definition}

\subsection{Known properties}
Definable families of sets in an o-minimal structure (such as those defined
above) have been studied and they satisfy important finiteness properties
similar to those of semi-algebraic families. We list here a couple
of properties which are important in the combinatorial study of 
arrangements.

\subsubsection{Finiteness of topological types}
\begin{theorem}[\cite{Dries,Michel2}]
\label{the:uniform}
Let ${\mathcal S}(\R)$ be an o-minimal structure over a real closed
field $\R$ and let $T \subset \R^{k+\ell}$
be a closed definable set. Then, the number of homeomorphism types
amongst $(T,\pi_1,\pi_2)$-sets is finite. 
\end{theorem}
\begin{remark}
\label{rem:uniform}
Note that, since the sum of the Betti number of any definable set is 
finite (since they are finitely triangulable \cite[Theorem 4.4]{Michel2})
Theorem \ref{the:uniform} implies that there  exists
a constant $C = C(T)$ (depending only on $T$) 
such that for any $(T,\pi_1,\pi_2)$-set $S$,
\[
\sum_{i=0}^{k} b_i(S) \leq C.
\]
\end{remark}
   
\subsubsection{Finiteness of VC dimension}
The notion of Vapnik-Chervonenkis dimension is important in many
applications in computational geometry (see \cite{Matousek}). 
We note here that 
$(T,\pi_1,\pi_2)$-families have finite Vapnik-Chervonenkis dimension,
for any fixed definable $T \subset \R^{k+\ell}$.
The following result is proved in \cite{Dries}.

We first recall the definition of the Vapnik-Chervonenkis dimension.
\begin{definition}[\cite{Matousek}]
\label{def:VC}
Let ${\mathcal F}$ be a set of subsets of an infinite set $X$. 
We say that a finite subset $A \subset X$ is shattered by ${\mathcal F}$
if each subset $B$ of $A$ can be expressed as $F_B \cap A$ for some
$F_B \in {\mathcal F}$. The VC-dimension of ${\mathcal F}$ is defined as
$$
\displaylines{
\sup_{A \subset X, |A| < \infty, A\; \mbox{is shattered by} \; 
{\mathcal F}} |A|.
}
$$
\end{definition}

\begin{theorem}[\cite{Dries}]
\label{the:VC}
Let $T$ be some definable subset of $\R^{k+\ell}$ 
in some o-minimal structure ${\mathcal S}(\R)$,
and
$\pi_1: \R^{k+\ell} \rightarrow \R^k,\pi_2: \R^{k+\ell} \rightarrow \R^\ell$ 
the two projections.
Then the VC-dimension of the family of
$(T,\pi_1,\pi_2)$-sets is finite.
\end{theorem}

\section{Main results}
\label{sec:results}
In this section we state our main results. As stated in the
Introduction, our goal is to study the 
combinatorial and topological complexity of sets defined in terms of
$n$ definable sets belonging to a fixed definable family in terms of
the parameter $n$. 
We show that the basic results
on combinatorial and topological complexity of arrangements continue
to hold in this setting.
Finally, as a sample application of our results
we extend a recent result of Alon et al.\cite{APPRS} 
on crossing
patterns of semi-algebraic sets to the o-minimal setting.

\begin{remark}
\label{rem:gp}
As remarked earlier, in many 
results on bounding  the combinatorial complexity of arrangements 
(of sets of constant description complexity) there is an assumption that
the sets be in general position \cite{AE98}. This is a rather strong
assumption and enables one to assume, for instance,
if the sets of the arrangements are hypersurfaces, 
that they intersect transversally and
this property usually plays a crucial role in the proof. 
In this paper we make no assumptions on general positions, 
nor on the objects of the arrangement themselves (apart from the fact that they
come from a fixed definable family).
The homological methods used in this paper make such assumptions unnecessary. 
\end{remark}

\subsection{Combinatorial and topological complexity of arrangements} 

\begin{theorem}
\label{the:betti}
Let ${\mathcal S}(\R)$ be an o-minimal structure over a real closed
field $\R$ and let $T \subset \R^{k+\ell}$
be a closed definable set. 
Then, there exists a constant $C = C(T) > 0$ depending only
on $T$, such that
for  any  $(T,\pi_1,\pi_2)$-family 
${\mathcal A} = \{S_1,\ldots,S_n\}$ 
of subsets of  $\R^k$ the following holds.
\begin{enumerate}
\item
For every $i, 0 \leq i \leq k$,
\[
\sum_{D \in {\mathcal C}({\mathcal A})} b_i(D) \leq C \cdot n^{k-i}.
\]
In particular, the combinatorial complexity of ${\mathcal A}$,
which is equal to 
\[
\sum_{D \in {\mathcal C}({\mathcal A})} b_0(D),
\]
is at most $C \cdot n^{k}.$   
\item
The topological complexity of any $m$ cells in the arrangement
${\mathcal A}$ is bounded by $m + C \cdot n^{k-1}$.
\end{enumerate}
\end{theorem}

Since dimension is a definable invariant (see \cite{Dries})
we can refine the notions of combinatorial  and topological complexity 
to arrangements restricted to a definable set of possibly smaller 
dimension than that of the ambient space as follows.

Let $V$ be a closed definable subset of $\R^k$ of dimension $k' \leq k$. 
For any $(T,\pi_1,\pi_2)$-family, ${\mathcal A} = \{S_1,\ldots,S_n\}$,
of  subsets of  $\R^k$,
and $I \subset \{1,\ldots,n\}$, we let ${\mathcal A}(I,V)$ 
denote the set
\begin{equation}
V \;\;\cap \;\;
\bigcap_{i \in I \subset [1\ldots n]} S_i  \;\; \cap 
\bigcap_{j \in [1\ldots n]\setminus I} 
\R^k \setminus S_j, 
\end{equation}
and we call a connected component of ${\mathcal A}(I,V)$ a
cell of the arrangement restricted to $V$.

Let ${\mathcal C}({\mathcal A},V)$ denote the set 
of all non-empty cells of the arrangement ${\mathcal A}$ restricted to $V$,
and we call
the cardinality of ${\mathcal C}({\mathcal A},V)$ the 
combinatorial complexity of the arrangement ${\mathcal A}$ restricted
to $V$. 
Similarly, we define the topological complexity
of an arrangement ${\mathcal A}$ restricted to $V$ to be the number
\[
\sum_{D \in {\mathcal C}({\mathcal A},V)} \sum_{i=0}^{k'} b_i(D).
\]

We have the following generalization of Theorem  \ref{the:betti}.
\begin{theorem}
\label{the:variety}
Let ${\mathcal S}(\R)$ be an o-minimal structure 
over a real closed field $\R$
and let 
$T \subset \R^{k+\ell}$,  $V \subset \R^k$
be closed definable sets with $\dim(V) = k'$. 
Then, there exists a constant $C = C(T,V) > 0$ depending only
on $T$ and $V$, such that
for any $(T,\pi_1,\pi_2)$-family, ${\mathcal A} = \{S_1,\ldots,S_n\}$,
of subsets of $\R^k$,
and for every $i, 0 \leq i \leq k'$,
\[
\sum_{D \in {\mathcal C}({\mathcal A},V)} b_i(D) \leq C \cdot n^{k'-i}.
\]  
In particular,
the combinatorial complexity of ${\mathcal A}$
restricted to $V$, which is equal to
$\sum_{D \in {\mathcal C}({\mathcal A},V)} b_0(D)$,
is bounded by $C \cdot n^{k'}$.
\end{theorem}

Now, let as before ${\mathcal S}(\R)$ be an o-minimal structure 
over a real closed field $\R$,
and let 
$T \subset \R^{k+\ell}$,  $V \subset \R^k$
be  closed definable sets with $\dim(V) = k'$.

\begin{theorem}
\label{the:A-sets}
Let ${\mathcal S}(\R)$ be an o-minimal structure 
over a real closed field $\R$,
and let 
$T \subset \R^{k+\ell}$,  $V \subset \R^k$
be closed definable sets with $\dim(V) = k'$. 
Then, there exists a constant $C = C(T,V) > 0$ such 
that for any  $(T,\pi_1,\pi_2)$-family, ${\mathcal A}$ with
$|\A| = n$,
and an ${\mathcal A}$-closed set $S_1 \subset \R^k$, and an
${\mathcal A}$-set $S_2 \subset \R^k$,
$$
\displaylines{
\sum_{i=0}^{k'} b_i(S_1 \cap V) \leq C \cdot n^{k'}\;\; \mbox{and,}\cr
\sum_{i=0}^{k'} b_i(S_2 \cap V) \leq C \cdot n^{2k'}.
}
$$
\end{theorem}

\subsection{Topological complexity of projections}
In Theorem \ref{the:A-sets} we obtained bounds on the topological complexity
of definable sets belonging to the Boolean algebra of sets generated by
any $(T,\pi_1,\pi_2)$-family of sets of cardinality $n$. We now consider
the images of such sets under linear projections. Such projections 
are closely related to the classical problem of quantifier elimination,
and play a very important role in semi-algebraic geometry.
In the case of semi-algebraic sets, there exist effective algorithms
for performing quantifier elimination, which enable one to compute
semi-algebraic descriptions of projections of semi-algebraic sets
in an efficient manner (see, for instance, \cite{BPRbook2}).
Notice however that unlike in the case of semi-algebraic sets, we do not
have effective algorithms for performing quantifier elimination over
a general o-minimal structure.

Using our theorem on quantitative cylindrical
definable cell decomposition (Theorem \ref{the:cdcd} below) it is possible
to give a doubly exponential bound (of the form
$C(T)\cdot n^{2(2^k -1)}$) on the sum of the Betti numbers of 
such projections. However, adapting a spectral sequence argument
from \cite{GVZ04}, we have the following singly exponential bound.
 
\begin{theorem}[(Topological complexity of projections)]
\label{the:projection}
Let 
${\mathcal S}(\R)$ 
be an o-minimal structure,
and let 
$T \subset \R^{k+\ell}$ 
be a definable, closed and bounded 
set. Let $k = k_1 + k_2$ and let 
$\pi_3: \R^k \rightarrow \R^{k_2}$ 
denote the projection map on the last $k_2$ co-ordinates.
Then, there exists a constant $C = C(T) > 0$ such 
that for any  $(T,\pi_1,\pi_2)$-family, ${\mathcal A}$, with
$|\A| = n$,
and an ${\mathcal A}$-closed set 
$S \subset \R^k$, 
$$
\displaylines{
\sum_{i=0}^{k_2} b_i(\pi_3(S)) \leq C \cdot n^{(k_1+1)k_2}.
}
$$
\end{theorem}

\subsection{Cylindrical definable cell decompositions}
In semi-algebraic geometry, {\em cylindrical algebraic decomposition}
is a very important method for obtaining a decomposition of an 
arbitrary semi-algebraic set into topological balls of various dimensions.
Once such a decomposition is computed, it can be refined to a semi-algebraic
triangulation, and various topological information about a given semi-algebraic
set (such as its Betti numbers) can be computed easily from 
such a triangulation.
Moreover, cylindrical algebraic decomposition can also be used for solving the 
quantifier elimination problem (see \cite{BPRbook2} for an exposition and
pointers to the large amount of literature on this subject).

The analogue of cylindrical algebraic decomposition over an o-minimal
structure is called {\em Cylindrical Definable Cell Decomposition}.
We first recall the definition of Cylindrical Definable Cell Decomposition
(henceforth called {\em cdcd}) following \cite{Michel2}.

\begin{definition}
\label{def:cdcd}
A cdcd of $\R^k$ is a finite partition of $\R^k$ into definable sets
$(C_i)_{i \in I}$ (called the cells of the cdcd) 
 satisfying the following properties.
\begin{enumerate}
\item
If $k=1$ then a cdcd of $\R$ is given by a finite set of points
$a_1 < \cdots < a_N$ and the cells of the cdcd are the singletons
$\{a_i\}$ as well as the open intervals,
$(\infty,a_1), (a_1,a_2),\ldots, (a_N,\infty)$.
\item
If $k > 1$, then a cdcd of $\R^k$ is given by a cdcd, $(C_i')_{i \in I'}$, 
of $\R^{k-1}$ and for each $i \in I'$,
a collection of cells, ${\mathcal C}_{i}$ defined by
\[
{\mathcal C}_{i} = \{ \phi_i(C_i' \times D_j) \;\mid\; j \in J_i\},
\]
where 
\[
\phi_i: C_i' \times \R \rightarrow \R^k
\]
is a definable homeomorphism satisfying $\pi \circ \phi = \pi$,
$(D_j)_{j \in J_i}$ is a cdcd of $\R$,
and
$\pi: \R^k \rightarrow \R^{k-1}$ is the projection map onto
the first $k-1$ coordinates.
The cdcd of $\R^k$ is then given by 
$$
\displaylines{
\bigcup_{i \in I'} {\mathcal C}_i.
}
$$
\end{enumerate}
Given a family of definable subsets $\A = \{S_1,\ldots,S_n\}$ of $\R^k$, 
we say that a cdcd is adapted to $\A$, if each $S_i$ is a union of cells
of the given cdcd.
\end{definition}

The fact that given any finite family $\A$ of definable subsets of $\R^k$, 
there exists a cdcd of $\R^k$ adapted to $\A$ is classical (see 
\cite{Michel2,Dries}). However, for the purposes of this paper we need a quantitative
version of this result.
Such quantitative versions are known in the semi-algebraic as well as
semi-Pfaffian categories (see, for example, \cite{BPRbook2,GVsurvey}), 
but is missing in the general o-minimal setting.

Given a $(T,\pi_1,\pi_2)$-family $\A$  
of cardinality  $n$,
we give a bound on the size of a cdcd of $\R^k$ adapted to this family
in terms of $n$, and furthermore show that cells of the cdcd come from
a definable family which depends only on $T$ (independent of $n$) and each
such cell can be defined only in terms of a constant number of elements
of $\A$. This latter property is essential in the combinatorial application 
described later in the paper. 

Since we will need to consider several  different projections, 
we adopt the following convention.
Given $m$ and $p$, $p\leq m$,  we will denote by 
$\pi_{m}^{\leq p}: \R^m \rightarrow \R^p$
(resp. $\pi_{m}^{>p}:\R^m \rightarrow \R^{m-p}$) 
the projection onto the first $p$ (resp. the last $m-p$) coordinates.

We prove the following theorem.

\begin{theorem}[(Quantitative cylindrical definable cell decomposition)]
\label{the:cdcd}
Let ${\mathcal S}(\R)$ be an o-minimal structure 
over a real closed field $\R$,
and let $T \subset \R^{k+\ell}$
be a closed definable set. 
Then, there exist constants $C_1,C_2 > 0$ depending only on $T$, 
and definable sets,
\[
\{T_\alpha\}_{\alpha \in I}, \;\; T_\alpha \subset \R^k \times \R^{2(2^k -1)\cdot\ell},
\]
depending only on $T$, 
with $|I| \leq C_1$, such that for any  $(T,\pi_1,\pi_2)$-family,
${\mathcal A} = \{S_1,\ldots,S_n\}$ with 
$S_i = T_{\y_i}, \y_i \in \R^{\ell}, 1 \leq i \leq n$, 
some sub-collection of the sets
$$
\displaylines{
\pi_{k+2(2^k -1)\cdot\ell}^{\leq k}
\left({\pi_{k+2(2^k -1)\cdot\ell}^{> k}}^{-1}(\y_{i_1},\ldots,\y_{i_{2(2^k -1)}}) 
\cap T_\alpha\right),\cr
\alpha \in I, \; 1 \leq i_1,\ldots,i_{2(2^k -1)} \leq n,
}
$$
form a cdcd of $\R^k$ compatible with $\A$. Moreover, the cdcd has
at most $C_2 \cdot n^{2(2^k -1)}$ cells.
\end{theorem}
 
The combinatorial complexity bound in 
Theorem \ref{the:cdcd} compares favorably with the combinatorial parts
of similar quantitative results on cylindrical decomposition of 
semi-algebraic sets (see for instance, Section 11.1 in \cite{BPRbook2}),
as well as sub-Pfaffian sets (see the main result in \cite{GV2}).
Moreover, since a doubly exponential dependence on $k$ is unavoidable 
(see \cite{DH}), the complexity bound in 
Theorem \ref{the:cdcd} is very close to the best possible. 
Notice also that it is possible to use Theorem \ref{the:cdcd} to give a 
doubly exponential bound on the Betti numbers of an ${\mathcal A}$-closed
set. However, we prove much better (singly exponential) bounds on the Betti
numbers of such sets (Theorems \ref{the:betti} and \ref{the:variety})
using different techniques.

\subsection{Application}
We end with an application (Theorem \ref{the:crossing} below)
which generalizes a Ramsey-type result due to
Alon et al. \cite{APPRS} from the class of semi-algebraic sets of
constant description complexity to $(T,\pi_1,\pi_2)$-families. 
One immediate consequence of Theorem \ref{the:crossing}
is that if we have two $(T,\pi_1,\pi_2)$-families,
${\mathcal A}$ and ${\mathcal B}$
of sufficiently large size, then one can always find a constant fraction,
${\mathcal A}' \subset {\mathcal A}$ , ${\mathcal B}' \subset {\mathcal B}$
of each, having the property that either every pair 
$(A,B) \in {\mathcal A}'\times {\mathcal B}'$ satisfy some
definable relation (for example, having a non-empty intersection) 
or no pair in ${\mathcal A}'\times {\mathcal B}'$ satisfy that relation.

More precisely, 

\begin{theorem}
\label{the:crossing}
Let ${\mathcal S}(\R)$ be an o-minimal structure 
over a real closed field $\R$,
and let $F$ be a closed definable subset of 
$\R^\ell \times \R^\ell$.
Then, there exists a constant $1 > \eps = \eps(F) > 0$, depending only
on $F$, such that for any set of $n$ points,
\[
{\mathcal F} = \{\y_1,\ldots,\y_n \in \R^{\ell} \}
\] 
there exists two subfamilies 
${\mathcal F}_1,{\mathcal F}_2 \subset {\mathcal F}$,
with $|{\mathcal F}_1|, |{\mathcal F}_2| \geq \eps n$
and either, 
\begin{enumerate}
\item
for all $\y_i \in {\mathcal F}_1$ and $\y_j \in {\mathcal F}_2$,
$(\y_i,\y_j) \in F$, or
\item
for no $\y_i \in {\mathcal F}_1$ and $\y_j \in {\mathcal F}_2$,
$(\y_i,\y_j) \in F$.
\end{enumerate}
\end{theorem}

An interesting application of Theorem \ref{the:crossing}
is the following.

\begin{corollary}
\label{cor:crossing}
Let ${\mathcal S}(\R)$ be an o-minimal structure 
over a real closed field $\R$, and let $T \subset \R^{k+\ell}$ be
a closed definable set.
Then, there exists a constant $1 > \eps = \eps(T) > 0$ depending only
on $T$, such that
for any  $(T,\pi_1,\pi_2)$-family, ${\mathcal A} = \{S_1,\ldots,S_n\}$,
there exists two subfamilies 
${\mathcal A}_1,{\mathcal A}_2 \subset {\mathcal A}$,
with $|{\mathcal A}_1|, |{\mathcal A}_2| \geq \eps n$,
and either, 
\begin{enumerate}
\item
for all $S_i \in {\mathcal A}_1$ and $S_j \in {\mathcal A}_2$,
$S_i \cap S_j \neq \emptyset$ or
\item
for all  $S_i \in {\mathcal A}_1$ and $S_j \in {\mathcal A}_2$,
$S_i \cap S_j = \emptyset$.
\end{enumerate}
\end{corollary}

\section{Proofs of the main results}
We first need a few preliminary results. 
\subsection{Finite unions of definable families} 
Suppose that $T_1,\ldots,T_m \subset \R^{k+\ell}$ are closed,
definable sets, $\pi_1: \R^{k+\ell} \rightarrow \R^{k}$ and
$\pi_2: \R^{k+\ell} \rightarrow \R^{\ell}$ the two projections.

We show that there exists a a certain closed definable 
subset $T' \subset \R^{k+\ell+m}$ 
depending only on $T_1,\ldots,T_m$, such that for any collection of 
$(T_i,\pi_1,\pi_2)$ families $\A_i$, $1 \leq i \leq m$,
the union, $\cup_{1\leq i \leq m} \A_i$, is a
$(T',\pi_1',\pi_2')$-family, 
where 
$\pi_1': \R^{k+m+\ell} \rightarrow \R^{k}$ and
$\pi_2': \R^{k+\ell+m} \rightarrow \R^{\ell+m}$ are the usual projections.

\begin{lemma}
\label{lem:unionoffamilies}
The family $\cup_{1 \leq i \leq m} \A_i$ is a $(T',\pi'_1,\pi'_2)$ family
where,
$$
T' = \bigcup_{i=1}^{m} T_i \times \{e_i\} \subset \R^{k+\ell + m},
$$
with $e_i$ the $i$-th standard basis vector in $\R^m$, and
$\pi_1': \R^{k+\ell+m} \rightarrow \R^k$ and
$\pi_2':   \R^{k+\ell+m} \rightarrow \R^{\ell+m},$  the
projections onto the first $k$ and the last $\ell+m$ coordinates 
respectively.
\end{lemma}

\begin{proof}
Obvious.
\end{proof}

\subsection{Hardt triviality for definable sets}
Our main technical tool will be
the following o-minimal version of Hardt's triviality theorem 
(see \cite{Dries,Michel2}).

Let $X \subset \R^k \times \R^\ell$ and $A \subset \R^\ell$ be
definable subsets  of $\R^k \times \R^\ell$ and $\R^\ell$ 
respectively,
and let $\pi: X \rightarrow \R^\ell$
denote the projection map.

We say that {\em $X$ is definably trivial over $A$} 
if there exists a definable
set $F$ and a definable homeomorphism 
\[
h: F \times A  \rightarrow X \cap \pi^{-1}(A), 
\]
such that the following diagram commutes:

\[
\xymatrix{
F \times A \ar[r]^h \ar[d]^{\pi_2} &  
X \cap \pi^{-1}(A) \ar[ld]^{\pi}\\
A&
},
\]
where $\pi_2: F \times A \rightarrow A$ is the projection onto
the second factor. We call $h$ 
{\em a definable trivialization of $X$ over $A$}.

If $Y$ is a definable subset of $X$, we say that the trivialization $h$ is
{\em compatible} with $Y$ if there is a definable subset $G$ of $F$ such
that $h(G \times A) = Y \cap \pi^{-1}(A)$. Clearly, the restriction of
$h$ to $G \times A$ is a trivialization of $Y$ over $A$.

\begin{theorem}[(Hardt's theorem for  definable families)]
\label{the:hardt}
Let $X \subset \R^k \times \R^\ell$ be a definable set and let
$Y_1,\ldots,Y_m$ be definable subsets of $X$. Then, there exists a 
finite partition of $\R^\ell$ into definable sets $C_1,\ldots,C_N$
such that $X$ is definably trivial over each $C_i$, and moreover
the trivializations over each $C_i$ are compatible with $Y_1,\ldots,Y_m$.
\end{theorem}

\begin{remark}
\label{rem:hardt}
Note that in particular it follows from Theorem \ref{the:hardt}, that
there are only a finite number of topological types amongst the
fibers of any definable map $f: X \rightarrow Y$ between definable
sets $X$ and $Y$ (see 
Remark \ref{rem:uniform}
).
\end{remark}

\subsection{Some notation}
For any definable set $X \subset \R^k$ we will denote by 
$X^c$  the complement of $X$, and by $\bar{X}$ the closure of $X$
in $\R^k$.
We also denote by $B_k(\x,r)$ (resp.
$\bar{B}_k(\x,r)$)
the open (resp. closed) ball in $\R^k$ of radius $r$ centered at $x$. 

For any closed definable subset $X \subset \R^k$,
we will denote by 
\[
d_X: \R^k \rightarrow \R, \;\; d_X(\x) = {\rm dist}(\x,X).
\]
Note that, it follows from the axioms in Definition \ref{def:o-minimal} 
that $d_X$ is a definable function (that is a function whose graph is 
a definable set).

Given closed definable sets $X \subset V \subset \R^k$,
and $\eps > 0$, 
we define the open tube of radius $\eps$ around $X$ in $V$ to 
be the definable set 
\[ 
\OT(X,V,\eps) = \{ \x \in V\;\mid\; d_X(\x) < \eps \}.
\] 

Similarly, we define the closed tube of radius $\eps$ around $X$ in $V$ to 
be the definable set 
\[ 
\CT(X,V,\eps) = \{ \x \in V\;\mid\; d_X(\x) \leq \eps \},
\]  
the boundary of the closed tube,
\[ 
\BT(X,V,\eps) = \{ \x \in V\;\mid\; d_X(\x) = \eps \},
\]  
and finally for $\eps_1 > \eps_2 > 0$ we define the open annulus
of radii $\eps_1,\eps_2$ around $X$ in $V$ to be the definable set
\[ 
\Ann(X,V,\eps_1,\eps_2) = \{ \x \in V \;\mid\; \eps_2 < d_X(\x) < \eps_1 \},
\]  
and the closed annulus
of radii $\eps_1,\eps_2$ around $X$ in $V$ to be the definable set
\[ 
\overline{\Ann}(X,V,\eps_1,\eps_2) = \{ \x \in V \;\mid\; \eps_2 \leq d_X(\x) \leq \eps_1 \}.
\]  

For brevity we will denote 
by $\OT(X,\R^k,\eps)$ (resp. $\CT(X,\R^k,\eps)$, $\BT(X,\R^k,\eps)$,
$\Ann(X,\R^k,\eps)$, $\overline{\Ann}(X,\R^k,\eps)$ ) by
$\OT(X,\eps)$ (resp.
 $\CT(X,\eps)$, $\BT(X,\eps)$, $\Ann(X,\eps)$,
$\overline{\Ann}(X,\eps)$).

\subsection{Replacing definable sets by closed and bounded 
ones maintaining homotopy type}

Let ${\mathcal A} = \{S_1,\ldots,S_n\}$ be a collection of closed, 
definable subsets of  $\R^k$ and let $V \subset \R^k$ be a closed,
and bounded definable set. 
In this section we adapt a construction due to Gabrielov and Vorobjov
\cite{GaV} for replacing any given $({\mathcal A},V)$-set 
by  a closed bounded $({\mathcal A}',V)$-set
(where ${\mathcal A}'$ is a new family of definable closely
related to ${\mathcal A}$)
such that the new set has the same homotopy type as the original one.

We denote by ${\rm In}({\mathcal A},V)$
the set,
\[
\{ I \subset [1\ldots n] \;\mid\; {\mathcal A}(I) \cap V \neq \emptyset\}.
\] 

Let,
$\varepsilon_{2n} \gg\varepsilon_{2n-1} \gg \cdots \gg \varepsilon_2 \gg 
\varepsilon_1 >0$  be sufficiently small.

For each $m, 0 \leq m \leq n$,
we denote by
${\rm In}_m({\mathcal A},V)$
the set $\{ I \in {\rm In}({\mathcal A},V) \;\mid\; |I| = m\}. $ 

Given $I \in {\rm In}_m({\mathcal A},V)$ denote by 
${\mathcal A}(I)^{cl}$ to be  the intersection of  $V$
with the closed definable set
\[
\bigcap_{i \in I} \CT(S_i,\eps_{2m}) \cap 
\bigcap_{i \in [1\ldots n]\setminus I} \overline{S_i^c}.
\]
and denote by
$\mathcal{A}(I)^{o}$ the intersection of  $V$ with 
the open definable set
\[
\bigcap_{i \in I} \OT(S_i,\eps_{2m-1}) \cap 
\bigcap_{i \in [1\ldots n]\setminus I} {S_i^c}.
\]

Notice that,
$$
\displaylines{
{\mathcal A}(I) \subset {\mathcal A}(I)^{cl}, \;\;\mbox{as well as}\cr
{\mathcal A}(I) \subset {\mathcal A}(I)^{o}.
}
$$

Let $X \subset  V$ be a $({\mathcal A},V)$-set
such that 
$\displaystyle{
X = \bigcup_{I \in \Sigma }{\mathcal A}(I) \cap V
}
$
with $\Sigma \subset {\rm In}({\mathcal A},V)$.
We denote
$\Sigma_m= \Sigma \cap {\rm In}_m({\mathcal A},V)$
and define a sequence of sets, 
$X^{m} \subset \R^k$, $0 \leq m \leq n$ inductively as follows.

\begin{itemize}
\item
Let $X^{0} = X.$
\item 
For 
$0 \leq m \leq n$,
we define 
$$
\displaylines{
X^{m+1} = 
\left(
X^{m} \cup  \bigcup_{I \in \Sigma_m} {\mathcal A}(I)^{cl}
\right)
\setminus 
\bigcup_{I \in {\rm In}_m({\mathcal A},V)
\setminus \Sigma_m}
{\mathcal A}(I)^{o}
}
$$
\end{itemize}

We denote by  $X'$ the set $ X^{n+1}$.

The following theorem is similar to Theorem 8.1 in
\cite{BPRbettione}. All the steps in the proof of Theorem 8.1 in 
\cite{BPRbettione} also remain valid in the o-minimal context. One needs 
to replace the references to Hardt's theorem for semi-algebraic mappings
by its o-minimal counterpart. Since repeating the entire proof with 
this minor modification would be tedious, we omit it from this paper.

\begin{theorem}
\label{the:GV}
The sets $X$ and $X'$ are  definably homotopy equivalent. 
\hide{
\begin{footnote}
{Very recently, after this paper was written,
Gabrielov and Vorobjov \cite{GV07} have given a much simpler
construction for replacing an arbitrary definable set $X$ by a closed and
bounded one, and if we use this new construction instead of the one described
above, we obtain a slightly improved bound in Theorem \ref{the:A-sets} 
(namely  
$C \cdot n^{k'}$ instead of $C \cdot n^{2k'}$). 
}
\end{footnote}
}
\end{theorem}
 
\begin{remark}
Very recently, after this paper was written,
Gabrielov and Vorobjov \cite{GV07} have given a much simpler
construction for replacing an arbitrary definable set $X$ by a closed and
bounded one, and if we use this new construction instead of the one described
above, we obtain a slightly improved bound in Theorem \ref{the:A-sets} 
(namely  
$C \cdot n^{k'}$ instead of $C \cdot n^{2k'}$). 
\end{remark}

\begin{remark}
\label{rem:GV}
Note that $X'$ is a 
$({\mathcal A}',V)$-closed set, 
where
\[
{\mathcal A}' = \bigcup_{i,j=1}^{n}\{S_i, 
\CT(S_i,\eps_{2j}), \OT(S_i,\eps_{2j-1})^c\}.
\]
If $\A$ is a $(T,\pi_1,\pi_2)$-family for some definable 
closed subset $T \subset \R^{k+\ell}$, then
by Lemma \ref{lem:unionoffamilies}, $\A'$ is 
a $(T',\pi_1',\pi_2')$-family for some definable
$T'$ depending only on $T$.
\end{remark}

\subsection{Mayer-Vietoris inequalities}
We will need a couple of inequalities which follows from 
the exactness of Mayer-Vietoris sequence. 
\begin{remark}
\label{rem:simplicial}
Note that for a closed and bounded definable set $X \subset \R^k$,
the homology groups $\HH_*(X)$ are isomorphic to the simplicial homology
groups of any definable triangulation of $X$ and in this case the proof
of the exactness of the Mayer-Vietoris sequence is purely combinatorial
in nature and presents no difficulties (even in the case when $\R$ is 
an arbitrary real closed field not necessarily equal to $\mathbb{R}$). 
The same remark also applies to arbitrary definable closed 
sets (not necessarily bounded), 
after intersecting the given sets
with  a large enough closed ball and using the conical structure at infinity
of definable sets.
\end{remark}

We first consider the case of two closed definable sets and
then generalize to the case of many such sets.

\begin{proposition}
\label{6:prop:MVc}
Let $S_1,S_2$ be two closed definable sets.
Then,

\begin{equation}
\label{6:eq:MV1}
b_i(S_1) + b_i(S_2) \leq b_i(S_1 \cup S_2) + b_i(S_1 \cap S_2),
\end{equation}
\begin{equation}
\label{6:eq:MV2}
b_i(S_1 \cup S_2) \leq b_i(S_1) + b_i(S_2) + b_{i-1}(S_1 \cap S_2),
\end{equation}
\begin{equation}
\label{6:eq:MV3}
b_i(S_1 \cap S_2) \leq b_i(S_1) + b_i(S_2) + b_{i+1}(S_1 \cup S_2).
\end{equation}
\end{proposition}

Let $S_1,\ldots, S_n \subset \R^k$ be closed definable sets,
contained in a closed bounded definable set $V$ of
dimension
$k'$.
For $1 \leq t \leq n$, we 
let 
\[
S_{\le t}= \bigcap_{1 \leq j \leq t} S_j, \;\mbox{and}\;
S^{\le t}= \bigcup_{1\leq j \leq t} S_j.
\]
Also,
for
$J \subset \{1,\ldots, n\}, \;J \neq \emptyset$, let
\[
S_J=\bigcap_{j \in J} S_j,\;\mbox{and}\; S^J=\bigcup_{j \in J} S_j. 
\]
Finally, let $S^\emptyset = V$.

We have the following proposition.
\begin{proposition}
\label{7:prop:prop1}
(a) For $0 \leq i \leq k'$,
\begin{equation}
b_i(S^{\le n}) \leq  \sum_{j=1}^{i+1}
\sum_{J \subset \{1,\ldots n\},\#(J)=j}
b_{i-j+1}(S_J).
\end{equation}

(b)
For $0\le i\le k',$
\begin{equation}
\label{7:eqn:prop1}
b_i(S_{\le n}) \leq
b_{k'}(S^\emptyset)+\sum_{j=1}^{k'-i} \;
\sum_{J \subset \{1,\ldots,n\}, \#(J)=j} \left(b_{i+j-1}(S^J) +
b_{k'}(S^\emptyset)\right).
\end{equation}
\end{proposition}

\begin{proof}
See \cite{BPR02}.
\end{proof}

\subsection{Proof of Theorem \ref{the:variety}}
We will use the following proposition in the proof of 
Theorem \ref{the:variety}.

\begin{proposition}
\label{prop:variety}
Let ${\mathcal A} = \{S_1,\ldots,S_n\}$ be a collection of closed 
definable subsets of  $\R^k$ and let $V \subset \R^k$ be a closed,
and bounded definable set. 
Then for all sufficiently small
$1 \gg \eps_1 \gg \eps_2 > 0$ the following holds.
For any connected component, $C$,  of ${\mathcal A}(I) \cap V$,
$I \subset [1\ldots n]$, there exists a
connected component, $D$,  of the definable set
\[ 
 \bigcap_{1 \leq i \leq n} \Ann(S_i,\eps_1,\eps_2)^c \cap V
\] 
such that $D$ is definably homotopy equivalent to $C$.
\end{proposition}

\begin{proof}
The proposition will follow from the following two observations  which are
consequences of Theorem \ref{the:hardt} 
(Hardt's theorem for o-minimal structures).

{\em Observation 1.}
It follows from Theorem \ref{the:hardt} that 
for all sufficiently small $\eps_1 > 0$ and 
for each connected component $C$ of 
${\mathcal A}(I) \cap V$, 
there exists a connected component $D'$ of 
\[
\bigcap_{i \in I} S_i \cap 
\bigcap_{j \in [1\ldots n]\setminus I}\OT(S_j,\eps_1)^c \cap V,
\]
definably homotopy equivalent to $C$.

{\em Observation 2.}
For all sufficiently small, $\eps_2$ with  $0 < \eps_2 \ll \eps_1$,
and for each connected component $D'$ of
\[
\bigcap_{i \in I} S_i \cap 
\bigcap_{j \in [1\ldots n]\setminus I}\OT(S_j,\eps_1)^c \cap V,
\]
there exists a connected component $D$ of
\[
W := \bigcap_{i \in I} \CT(S_i,\eps_2) \cap 
\bigcap_{j \in [1\ldots n]\setminus I}\OT(S_j,\eps_1)^c \cap V,
\]
definably
homotopy equivalent to $D'$.

Now notice that $D$ is connected and contained in the set
\[ 
\bigcap_{1 \leq i \leq n} \Ann(S_i,\eps_1,\eps_2)^c \cap V.
\]
Let $D''$ be the connected component of 
\[ 
\bigcap_{1 \leq i \leq n} \Ann(S_i,\eps_1,\eps_2)^c \cap V 
\]
containing $D$. We claim that $D = D''$, which will prove the
proposition.
 
Suppose $D'' \setminus D \neq \emptyset$. 
Then, $D'' \setminus W \neq \emptyset$, since otherwise $D'' \subset W$,
which would imply that $D'' = D$, since $D''$ is connected and $D\subset D''$ 
is a connected component of $W$.
Let $\x \in D'' \setminus W$
and $\y$ any point in $D$. 
Since $\x \notin W$, either
\begin{enumerate}
\item
there exists $i \in I$ such that $\x \in \OT(S_i,\eps_1)^c$ or

\item 
there exists $i \in [1\ldots n]\setminus I$ such that 
$\x \in \CT(S_i,\eps_2)$.
\end{enumerate}

Let $\gamma: [0,1] \rightarrow D''$ be a definable path with
$\gamma(0) = \x, \gamma(1) = \y$.
and let $d_i: D'' \rightarrow \R$ be the definable continuous function,
$d_i(\z) = {\rm dist}(\z,S_i).$

Then, in the first case, $d_i(\x) =  d_i(\gamma(0)) \geq \eps_1$ and 
$d_i(\y) = d_i(\gamma(1)) <  \eps_2$, implying that there exists
$t \in (0,1)$ with $\eps_2 <  d_i(\gamma(t)) <  \eps_1$ implying that
$d_i(\gamma(t)) \not\in \Ann(S_i,\eps_1,\eps_2)^c$ and hence not in 
$D''$ (a contradiction).
In the second case,
$d_i(\x) =  d_i(\gamma(0)) <  \eps_2$ and 
$d_i(\y) = d_i(\gamma(1)) \geq   \eps_1$, implying that there exists
$t \in (0,1)$ with $\eps_2 <  d_i(\gamma(t)) < \eps_1$ again implying that
$d_i(\gamma(t)) \not\in \Ann(S_i,\eps_1,\eps_2)^c$ and hence not in 
$D''$ (a contradiction).
\end{proof}

We are now in a position to prove Theorem \ref{the:variety}.

\begin{proof}[
of Theorem \ref{the:variety}]
For $1 \leq i \leq n$, let $\y_i \in \R^{\ell}$ such that
\[
S_i =  T_{\y_i},
\]
and let 
\[
A_i(\eps_1,\eps_2) = \Ann(S_i,\eps_1,\eps_2)^c \cap V.
\]

Applying Proposition \ref{7:prop:prop1} we have for $0 \leq i \leq k'$,
\begin{equation}
\label{eqn:the:variety}
b_i(\bigcap_{j=1}^{n}A_j(\eps_1,\eps_2)) \leq
b_{k'}(V)+\sum_{j=1}^{k'-i} \;
\sum_{J \subset \{1,\ldots,n\}, \#(J)=j} \left(b_{i+j-1}(A^J(\eps_1,\eps_2)) +
b_{k'}(V)\right),
\end{equation}
where $A^J(\eps_1,\eps_2) = \cup_{j \in J} A_j(\eps_1,\eps_2)$.

Notice that each $\Ann(S_i,\eps_1,\eps_2)^c, 1\leq i \leq n$, is a
$(\Ann(T,\eps_1,\eps_2)^c,\pi_1,\pi_2)$-set and moreover,
\[
\Ann(S_i,\eps_1,\eps_2)^c = \pi_1(\pi_2^{-1}(\y_i)\cap 
\Ann(T,\eps_1,\eps_2)^c) ;\;1 \leq i \leq n.
\] 

For $J \subset [1\ldots n]$
with $|J| \leq k'$, 
we will denote 
\[
S^J(\eps_1,\eps_2) = \bigcup_{j \in J} \Ann(S_j,\eps_1,\eps_2)^c.
\]

Consider the definable set
\[
B_J(\eps_1,\eps_2) = \prod_{j \in J} \Ann(T,\eps_1,\eps_2) \cap \Delta,
\]
where 
$\Delta \subset \R^{|J|(k+\ell)}$ is the definable (in fact, semi-algebraic) 
set defined by 
\[
\Delta = \{(\x,\z_1,\x,\z_2,\ldots,\x,\z_{|J|}) \; \mid\; \x \in \R^{k}, 
\z_1,\ldots,\z_{|J|} \in \R^{\ell}\}.
\]

The projection map $\pi_2$ induces a projection map,
\[
\prod_{j \in J} \pi_2: \R^{|J|(k+\ell)} \rightarrow \prod_{j \in J} \R^{\ell}.
\]

We also have the natural projection 
\[
\pi_1: B_J(\eps_1,\eps_2)  \rightarrow \R^{k}.
\]
 
\[
\xymatrix{
&  B_J(\eps_1,\eps_2) \ar[ld]^{\pi_1}
\ar[rd]^{ \prod_{j \in J} \pi_2}  & \\
\R^{k} & & \R^{|J|\ell}
}
\]

It is now easy to see that for each, 
$J = \{i_1,\ldots,i_{|J|}\}$,
$S^J(\eps_1,\eps_2)^c$ is homeomorphic to
$(\prod_{j\in J} \pi_2)^{-1}(\y_{i_1},\ldots,\y_{i_{|J|}}) 
\cap B_J(\eps_1,\eps_2)$ 
via the projection $\pi_1$.

Using Remark \ref{rem:hardt} we can conclude
there exists an upper bound depending only on $T$ (and independent of
$\y_1,\ldots,\y_n$ as well as $\eps_1,\eps_2$), on the number of 
of topological types 
amongst the pairs  
\[
\left( \R^k, \pi_1\left((\prod_{j \in J}\pi_2)^{-1}(\y_{i_1},\ldots,\y_{i_{|J|}}) \cap 
B_J(\eps_1,\eps_2)\right)\right),
\]
and hence amongst the pairs  $(\R^k, S^J(\eps_1,\eps_2)^c)$ as well.
This implies that  
there are only a finite number (depending on $T$) of topological types 
amongst $S^J(\eps_1,\eps_2)$.
Restricting all the sets to $V$ in the above argument, 
we obtain that there are only finitely many
(depending on $T$ and $V$) of topological types amongst the sets
$A^J(\eps_1,\eps_2) =  S^J(\eps_1,\eps_2) \cap V$.

Thus, there exists a constant $C(T,V)$ such that
$$\displaylines{
C(T,V) =  \max_{J \subset \{1,\ldots,n\}, |J| \leq k', 
0 \leq i+j \leq k'} \left(b_{i+j-1}(A^J(\eps_1,\eps_2)) +
b_{k'}(V)\right) + b_{k'}(V).
}
$$
It now follows from inequality \ref{eqn:the:variety} and Proposition
\ref{prop:variety} that,
\[
\sum_{D \in {\mathcal C}({\mathcal A},V)} b_i(D) \leq C \cdot n^{k'-i}.
\]   
\end{proof}

We now prove Theorem \ref{the:A-sets}.

The proof of Theorem \ref{the:A-sets} will follow from the following
proposition.
For the sake of greater clarity,
and since it does not affect in any way the proof of Theorem \ref{the:A-sets}, 
we choose to be slightly less precise in the next proposition, 
and prove a bound on the sum of the Betti numbers of $S$
(rather than  prove separate bounds on each individual
Betti number).
Recall from before that for any  definable set $X \subset \R^k$,
we denote by $b(X)$ the sum  $\sum_{i\geq 0} b_i(X)$.

\begin{proposition}
\label{prop:variety2}
Let ${\mathcal A} = \{S_1,\ldots,S_n\}$ be a collection of closed 
definable subsets of  $\R^k$ and let $V \subset \R^k$ be a closed,
and bounded definable set and let $S$ be an $({\mathcal A},V)$-closed set.
Then, for all sufficiently small
$1 \gg \eps_1 \gg \eps_2 \cdots \gg \eps_n > 0$,
\[
b(S) \leq
 \sum_{D \in {\mathcal C}({\mathcal B},V)} b(D),
\] 
where 
\[
{\mathcal B} = \bigcup_{i=1}^{n}\{ S_i, \BT(S_i,\eps_i), 
\OT(S_i,2\eps_i)^c
\}.
\]
\end{proposition}

\begin{proof}[
of Proposition \ref{prop:variety2}]
We define ${\mathcal A}_{>i}=\{S_{i+1},\ldots, S_n\}$
and
\[
{\mathcal B}_i=  \{ S_i, \BT(S_i,\eps_i),\OT(S_i,2\eps_i)^c \},
\]
and 
\[
{\mathcal B}_{\le i}= \{X \mid
X =\bigcap_{j=1,\ldots,i} X_j, X_j \in {\mathcal B}_j\}.
\]

The proof of the proposition will follow from the following proposition.

\begin{proposition}
\label{7:prop:closed}
For every $({\mathcal A},V)$-closed set $S$,
$$
b(S)  \leq
\sum_{X \in {\mathcal B}_{\le s}, X \cap V \subset S } b(X \cap V).
$$
\end{proposition}

The main ingredient of the 
proof of the proposition is the following lemma.

\begin{lemma}
\label{7:lem:closed}
For every $({\mathcal A},V)$-closed set $S$,
 and every $X \in {\mathcal B}_{\le i}$,
$$
b(S \cap X)  \leq
\sum_{Y \in {\mathcal B}_{i+1} } b(S \cap X \cap Y).
$$
\end{lemma}

\begin{proof}[
of Lemma \ref{7:lem:closed}]
Consider the sets
$$\displaylines{
T_1=S \cap X \cap \OT(S_{i+1},\eps_{i+1})^c,\cr
T_2=S \cap X \cap  \CT(S_{i+1},3\eps_{i+1}).
}
$$

Clearly, $S \cap  X = T_1 \cup T_2$.

Using Proposition \ref{6:prop:MVc}, we have that,
\[
b(S \cap X)  \leq b(T_1) + b(T_2) +
b(T_1 \cap  T_2).
\]

Now, since 
\[
T_1 \cap  T_2 = S \cap  X \cap  \overline{\Ann}(S_{i+1},3\eps_{i+1},\eps_{i+1})),
\]
we have that,
\[
b(T_1 \cap T_2)=b(S \cap X \cap  \overline{\Ann}(S_{i+1},3\eps_{i+1},\eps_{i+1})).
\]

It is now easy to verify using Theorem \ref{the:hardt} that,
$$
\displaylines{
T_1 \sim S \cap X \cap \OT(S_{i+1},2\eps_{i+1})^c,\cr
T_2 \sim S \cap X \cap S_{i+1},\cr
T_1 \cap T_2 \sim S \cap X \cap \BT(S_{i+1},2\eps_{i+1}),
}
$$
where $\sim$ denotes definable homotopy equivalence.

Finally,
$$
b(S \cap X)  \leq
\sum_{Y \in {\mathcal B}_{i+1}}  b(S\cap  X \cap Y).
$$
\end{proof}

\begin{proof}[
of Proposition \ref{7:prop:closed}]
Starting from the set  $S$ apply Lemma
\ref{7:lem:closed} with $X$ the empty set. Now,
repeatedly apply Lemma \ref{7:lem:closed} to the terms
appearing on the right-hand side of the inequality obtained,
noting that for any
$Y \in {\mathcal B}_{\le s},$
either $S \cap X= X$, and thus $ X \subset S$,
or $S \cap  X=\emptyset$.
\end{proof}

The proof of Proposition \ref{prop:variety2} now follows from 
Proposition \ref{7:prop:closed}.
\end{proof} 

\subsection{Proof of Theorem \ref{the:A-sets}}

\begin{proof}[
of Theorem \ref{the:A-sets}]
Follows directly from Theorem \ref{the:GV}, Theorem \ref{the:variety} and 
Proposition \ref{prop:variety2}.
\end{proof}

\subsection{Proof of Theorem \ref{the:projection}}
The proof of Theorem~\ref{the:projection} relies on the bounds in
Theorem~\ref{the:betti}, and on the following theorem which is
adapted to the o-minimal setting from \cite{GVZ04}.

\begin{theorem}
\label{the:ss}
Let $X$ and $Y$ be two closed, definable sets and $f: X \to Y$ a
definable continuous surjection 
which is closed (i.e. $f$ takes closed sets to closed sets).
Then for any integer $q$, we have

\begin{equation}\label{eqn:ss}
b_q(Y) \leq \sum_{i+j=q} b_j(W^i_f(X)),
\end{equation}
where $W^i_f(X)$ denotes the $(i+1)$-fold fibered product of $X$ over $f$:
\[
W^i_f(X) = \{(\x_0,\ldots,\x_i) \in X^{i+1}
\mid f(\x_0) = \cdots = f(\x_i)\}.
\]
\end{theorem}

\begin{remark}(Regarding the proof of Theorem  \ref{the:ss}.)
Theorem \ref{the:ss} was proved in \cite{GVZ04} in the semi-algebraic
and semi-Pfaffian setting and follows from the existence of a spectral
sequence $E_r^{i,j}$ converging to $H^*(Y)$ and such that $E_1^{i,j}
\cong H^j(W^i_f(X))$. Thus, the extension of this theorem to general
o-minimal structures over arbitrary real closed fields $\R$ (not necessarily
equal to $\mathbb{R}$) requires some remarks. 
The existence of the spectral sequence $E_r^{i,j}$ is a consequence of the
$p$-connectivity of the $(p+1)$-fold join of any simplicial complex $K$, 
and the Vietoris-Begle theorem.
The proof of the 
$p$-connectivity of the $(p+1)$-fold join of any simplicial complex $K$
is combinatorial in nature (see, for instance, 
\cite[Proposition 4.4.3]{Matousek_book2}), and thus 
presents no additional difficulties over general o-minimal structures.
A purely combinatorial proof of the Vietoris-Begle theorem 
is also known \cite[Theorem 2]{BWW06} (see also \cite[Corollary 2.6]{GV07}).
Since the rest of the argument is combinatorial in nature, 
it extends without 
difficulty to closed maps in  arbitrary o-minimal structures after
choosing appropriate triangulations.
Finally, since in any spectral sequence, the dimensions of
the terms $E_r^{i,j}$ are non-increasing when $i$ and $j$ are fixed and
$r$ increases, we obtain:
\[
b_n(Y)=
\sum_{i+j=n} \dim \left(E_{\infty}^{i,j}\right)
\leq \sum_{i+j=n} \dim \left(E_1^{i,j}\right),
\]
yielding inequality~(\ref{eqn:ss}).
\end{remark}

\begin{proof}[
of Theorem \ref{the:projection}]
Notice that for each $p, 0 \leq p \leq k_2$, and any ${\mathcal A}$-closed
set $S \subset \R^{k_1+k_2}$,
$W^p_{\pi_3}(S) \subset \R^{(p+1)k_1+ k_2}$ 
is an ${\mathcal A}^p$-closed set where,
$$
\displaylines{
{\mathcal A}^p = \bigcup_{j=0}^{p}{\mathcal A}^{p,j}, \cr
{\mathcal A}^{p,j} = \bigcup_{i=1}^{n}\{S_i^{p,j}\},
}
$$
where $S_i^{p,j} \subset \R^{(p+1)k_1 + k_2}$ is defined by
$$
\displaylines{
S_i^{p,j} = 
\{(\x_0,\ldots,\x_p,\y) \; \mid\; \x_j \in \R^{k_1}, \y \in \R^{k_2},
(\x_j,\y) \in S_i\}.  
}
$$
Also, note that ${\mathcal A}^{p,j}$ is a $(T^{p,j},\pi_1^p,\pi_2^p)$ 
family,
where 
$$
\displaylines{
T^{p,j} = \{(\x_0,\ldots,\x_p,\y,\z) \; \mid\; \x_j \in \R^{k_1}, 
\y \in \R^{k_2},
\z \in \R^\ell, (\x_j,\y,\z) \in T, \;\cr
\mbox{for some}\; j, 0 \leq j \leq p\}.  
}
$$
and
$\pi_1^p: \R^{(p+1)k_1 + k_2 + \ell} \rightarrow \R^{(p+1)k_1+k_2}$,
and 
$\pi_2^p: \R^{(p+1)k_1 + k_2 + \ell} \rightarrow \R^{\ell}$
are the appropriate projections. 
Since each $T^{p,j}$ is determined by $T$, we have  using
Lemma \ref{lem:unionoffamilies} that
$\A^p$ is a $(T',\pi_1',\pi_2')$-family for some definable 
$T'$ determined by $T$.
Note that, $W^p_{\pi_3}(S) \subset \R^{(p+1)k_1+ k_2}$ is a 
${\mathcal A}^p$-closed set and $\#{\mathcal A}^p = (p+1)n$.
Applying Theorem \ref{the:betti} we get, for each $p$ and $j$, 
$0 \leq p,j < k_2$,
$$
\displaylines{
b_j(W^p_{\pi_3}(S)) \leq C_1(T) \cdot n^{(p+1)k_1 + k_2}
}
$$
The theorem now follows from Theorem \ref{the:ss}, since
for each $q, 0 \leq q < k_2$,
\[
b_q(\pi_3(S)) \leq \sum_{i+j=q} b_j(W^i_{\pi_3}(S)) \leq C_2(T)\cdot
n^{(q+1)k_1 + k_2} \leq C(T)\cdot n^{(k_1+1)k_2}.
\]
\end{proof}

\subsection{Proof of Theorem \ref{the:cdcd}}
The proof of Theorem \ref{the:cdcd} will follow from the following lemma
(which corresponds to the first projection step in the 
more familiar cylindrical algebraic decomposition algorithm
for semi-algebraic sets (see, for instance, \cite{BPRbook2})).

\begin{lemma}
\label{lem:cdcd}
Let ${\mathcal S}(\R)$ be an o-minimal structure over a real closed field $\R$,
and let $T \subset \R^{k+\ell}$ be a closed definable set.  
Then, there 
exists definable sets $T_1,\ldots,T_N \subset \R^{k-1+2\ell}$ satisfying
the following.
For  each $i, 1 \leq i \leq N$, and $\y,\y' \in \R^{\ell}$, 
let 
$$
\displaylines{
B_i(\y,\y') = \pi_{k-1+2\ell}^{\leq k-1}
\left(
{\pi_{k-1+2\ell}^{> k-1}}^{-1}
(\y_1,\y_2) \cap T_i
\right).
}
$$
The projection
$\pi_{k}^{>1}: \R^{k} \rightarrow \R^{k-1}$ restricted to the sets
$T_{\y} \cup T_{\y'}$ 
is definably  trivial over $B_i(\y,\y')$ 
and the trivialization
is compatible with $T_{\y}$ and $T_{\y'}$.
\end{lemma}

\begin{proof}
Let 
$$
\displaylines{
V_0 = \{ (\x,\y,\y') \;\mid\; (\x,\y) \in T \;\mbox{or} \; (\x,\y') 
\in T \}, \cr
V_1 = \{ (\x,\y,\y') \;\mid\; (\x,\y) \in T \}, \cr
V_2 = \{ (\x,\y,\y') \;\mid\; (\x,\y') \in T \}.
}
$$
Note that $V_0 \subset \R^{k+2\ell}$ and $V_1,V_2 \subset V_0$ and
$V_0,V_1,V_2$ are all definable and determined by $T$. 
Applying Hardt's triviality theorem to the
sets $V_0,V_1,V_2$ and the projection map $\pi_{k+2\ell}^{> 1}$, we get a definable
partition of $\R^{k-1+2\ell}$ into definable sets $T_1,\ldots,T_N$,
such that $\pi_{k+2\ell}^{> 1}|_{V_0}$ can be trivialized over each $T_i$ and the
trivializations respects the subsets $V_1,V_2$. It is now easy to check that
the sets $T_i$ have the required properties.
\end{proof}

\begin{proof}[
of Theorem \ref{the:cdcd}]
We will use induction on $k$. 

The base case is when $k=1$ and the theorem is clearly true in this case.

Now suppose by induction hypothesis that the theorem is true for $k-1$.
We first apply Lemma \ref{lem:cdcd} to obtain 
definable sets $T_1,\ldots,T_N \subset \R^{k-1+2\ell}$ satisfying
the following conditions. 

For each $i, 1 \leq i \leq N$, and $\y,\y' \in \R^{\ell}$, 
the projection
$\pi_{k}^{>1}: \R^{k} \rightarrow \R^{k-1}$ restricted to the sets
$T_{\y}\cup T_{\y'}$ is definably trivial over $B_i(\y,\y')$ 
and the trivialization is compatible with $T_{\y}$ and $T_{\y'}$,
where
$$
\displaylines{
B_i(\y,\y') = 
\pi_{k-1+2\ell}^{\leq k-1}\left
({\pi_{k-1+2\ell}^{> k-1}}^{-1}
(\y,\y') \cap T_i
\right),
}
$$

Now let $T' = \cup_{1 \leq i \leq N} B_i \times \{e_i\}$ where $e_i$ is
that $i$-th standard basis vector in $\R^N$. Note that
$T' \subset \R^{k-1+ 2\ell + N}$. 

Applying the induction hypothesis to the triple
$$
\left(
T' \subset \R^{k-1 + 2(2^{k-1} -1)\cdot (2\ell+N)},
\pi_{k-1 + 2(2^{k-1} -1)\cdot (2\ell+N)}^{\leq k-1},
\pi_{k-1 + 2(2^{k-1} -1)\cdot (2\ell+N)}^{> k-1}
\right)
$$
we obtain definable sets,
\[
\{T_j'\}_{j \in J}, \;\; T_j' \subset \R^{k-1} \times \R^{2(2^{k-1}-1)\cdot(2\ell + N)},
\]
depending only on $T$ having the property that,
for any
$\y_1,\ldots,\y_n, \in \R^{\ell}$
and $\a = (\a_1,\ldots,\a_{2(2^{k-1} -1)}) \in \R^{2(2^{k-1} -1)\cdot N}$
where each $\a_i$ is a standard basis vector in $\R^N$,
some sub-collection of the sets
$$
\displaylines{
\pi_{k-1 + 2(2^{k-1} -1)\cdot (2\ell+N)}^{\leq k-1}
\left({\pi_{k-1 + 2(2^{k-1} -1)\cdot (2\ell+N)}^{> k-1}}^{-1}(\y_{i_1},
\ldots,\y_{i_{2^2(2^k -1)}},\a) 
\cap T_i\right),
}
$$
form a cdcd of $\R^k$ compatible with the family
$$
\displaylines{
\bigcup_{1 \leq i,j  \leq n}\bigcup_{1 \leq h \leq N}
\{B_h(\y_i,\y_j)\}.
}
$$

For $\x \in \R^{k-1}$, and  $\y\in \R^{\ell}$,
let
$$
\displaylines{
S(\x,\y) = \{x \in \R \;\mid\; (x,\x,\y) \in T \}.
}
$$

Now,  for $\x \in \R^{k-1}, \y,\y'\in \R^{\ell}$,
$S(\x,\y),S(\x,\y') \subset \R$ and each of them is a union
of a finite number of open intervals and points. 
The sets $S(\x,\y),S(\x,\y')$ induce  a partition of $\R$ 
into pairwise disjoint subsets, 
\[
V_1(\x,\y,\y'),V_2(\x,\y,\y'),\ldots,
\] 
where for $i \geq 0$, each
$V_{2i+1}(\x,\y,\y')$ is a maximal open interval contained
in one of 
$$
\displaylines{
S(\x,\y)\cap S(\x,\y'), S(\x,\y)^c \cap S(\x,\y'), \cr
S(\x,\y) \cap S(\x,\y')^c, S(\x,\y)^c \cap S(\x,\y')^c,
}
$$
and  
$V_{2i}(\x,\y,\y')$ is the right  end-point of the interval 
$V_{2i-1}(\x,\y,\y')$. 
We let
${\mathcal V}(\x,\y,\y')$ denote the ordered sequence,
\[
\left(V_1(\x,\y,\y'),V_2(\x,\y,\y'),\ldots,V_M(\x,\y,\y')\right),
\]
where $M$ is a uniform upper bound on $|{\mathcal V}|$
depending on $T$, and with the
understanding that $V_i(\x,\y,\y')$ can be empty for all $i \geq i_0$
for some $0 \leq i_0 \leq M$.
It is clear that the sets,
\[
V_i = \{(V_i(\x,\y,\y'),\x,\y,\y') \;\mid\; \x \in \R^{k-1},\y,\y'
\in \R^{\ell}
\}
\] 
are definable and depend only on $T$.    

For each $T_j' \subset \R^{k-1} \times \R^{2(2^{k-1}-1)\cdot(2\ell + N)},
 j \in  J$,
$1 \leq h \leq M$, 
and 
$\a = (\a_1,\ldots,\a_{2^{k} -2}),
$ 
where each $\a_i$ is a standard basis vector in $\R^N$,
let
$$
\displaylines{
T_{j,h,\a}' = \{(V_h(\x,\y_{2^{k+1} -3},\y_{2^{k+1} -2}),
\x,\y_1,\ldots,\y_{(2(2^k - 1)}) \; \mid\cr
\; (\x,\y_1,\ldots,\y_{2^2(2^{k-1}-1)},\a) \in T_j'
\}.
}
$$
Let $\{T_i\}_{i \in I}$ be the collection of all possible
$T_{j,h,\a}'$. It is now easy to verify that the family of sets
$\{T_i\}_{i \in I}$ satisfies the conditions of the theorem. 
\end{proof}

\subsection{Proof of Theorem \ref{the:crossing}}
The proof is very similar to the second proof of Theorem 1.1 in \cite{APPRS}.
However, instead of using vertical decomposition as in \cite{APPRS},
we use the cylindrical definable cell decomposition 
given by Theorem \ref{the:cdcd}. We repeat it here for the reader's
convenience.

\begin{proof}[
of Theorem \ref{the:crossing}]
For each $i, 1 \leq i \leq n$, let 
\[
A_i = \pi_{2\ell}^{\leq \ell}({\pi_{2\ell}^{> \ell}}^{-1}(\y_i) \cap F),
\]
and ${\mathcal G} = \{A_i \;\mid\; 1 \leq i \leq n \}$.
Note that ${\mathcal G}$ is a 
$(R,\pi_{2\ell}^{\leq \ell}, \pi_{2\ell}^{> \ell})$-family.

We now use  the Clarkson-Shor random sampling technique \cite{Matousek} 
(using Theorem \ref{the:cdcd} instead of vertical decomposition 
as in \cite{APPRS}).
Applying Theorem \ref{the:cdcd} to some sub-family 
${\mathcal G}_0 \subset {\mathcal G}$ of
cardinality $r$,
we get a decomposition of $\R^\ell$ into at most 
$C r^{2(2^{\ell}-1)} = r^{O(1)}$ 
definable cells, each of them defined by at most 
$2(2^{\ell}-1)= O(1)$ of the $\y_i$'s.
This decomposition satisfies the necessary properties for the existence of
$1/r$-cuttings of size $r^{O(1)}$ \cite[pp. 163]{Matousek}.

More precisely,
let $\tau$ be a cell of the cdcd of ${\mathcal G}_0$ and let 
$G \in {\mathcal G}$. We say that $G$ crosses $\tau$ if 
$G \cap \tau \neq \emptyset$ and $\tau \not\subset G$.
The well-known Cutting Lemma (see \cite[Chapter 6, Section 5]{Matousek})
now ensures that we can 
choose ${\mathcal G}_0$ such that each cell of the cdcd of ${\mathcal G}_0$
is crossed by no more than $\frac{c_1 n \log r}{r}$ elements of 
${\mathcal G}$, where $c_1$ is a constant depending only on $F$.

For each cell $\tau$ of the cdcd of ${\mathcal G}_0$, let
${\mathcal G}_\tau$ denote  the set of elements of ${\mathcal G}$ which
cross $\tau$ and let ${\mathcal F}_\tau = {\mathcal F} \cap \tau$.

Since the total number of cells in the cdcd of ${\mathcal G}_0$ is bounded by
$r^{O(1)}$,
there must exist a cell $\tau$ such that,
\[
|{\mathcal F}_\tau| \geq \frac{n}{r^{O(1)}}.
\]
Now, every element of ${\mathcal G} \setminus {\mathcal G}_\tau$ either
fully contains $\tau$ or is disjoint from it.

Setting $\alpha = \frac{1}{r^{O(1)}}$ and 
$\beta = \frac{1}{2}(1 - \frac{c_1 \log r}{r})$
we have that there exists a set ${\mathcal F}' = {\mathcal F}_\tau$
of cardinality at least $\alpha n$, and a subset 
${\mathcal G}'$ of cardinality at least
 $\beta n$ such that either each element of ${\mathcal F}'$
is contained in every element of ${\mathcal G}'$, or no
element of ${\mathcal F}'$
is contained in any  element of ${\mathcal G}'$.

The proof is complete by taking ${\mathcal F}_1 = {\mathcal F}'$,
and ${\mathcal F}_2 = \{\y_i\; \mid \; A_i \in {\mathcal G}'\}$ and
choosing $r$ so as to maximize
$\eps = \min(\alpha,\beta)$.
\end{proof}

\begin{proof}[
of Corollary  \ref{cor:crossing}]
For $1 \leq i \leq n$, let $\y_i \in \R^{\ell}$ be such that 
$S_i = T_{\y_i}$.
Let $F \subset \R^{\ell} \times \R^{\ell}$ be the closed 
definable set defined by
\[ 
F = \{(\z_1,\z_2) \;\mid\; \z_1,\z_2 \in \R^{\ell},\;\;
T_{\z_1} \cap T_{\z_2} 
\neq \emptyset\}.
\]

Clearly, $F$ is completely determined by $T$. Now apply 
Theorem \ref{the:crossing}.
\end{proof}

\section{Conclusion and open problems}
In this paper we have proved bounds on the combinatorial and topological 
complexities of arrangements of sets belonging to some fixed definable
family in an o-minimal structure, in terms of the number of sets in the
arrangement. These results generalize known results in the case when
the sets in the arrangements are semi-algebraic sets and of 
constant description complexity.
We also extended a Ramsey-type theorem due to Alon et al. \cite{APPRS}, 
originally proved for semi-algebraic sets of fixed description complexity 
to the more general setting of o-minimal geometry.

There are many other sophisticated results on the combinatorial complexity
of sub-structures of arrangements which have been proved in the semi-algebraic
case. Usually there are some extra assumptions about general position in these
results. For instance, it was shown in \cite{Basu3} that the complexity
of a single cell in an arrangement of $n$ semi-algebraic hyper-surface 
patches in $\R^k$, which are in general position and have  
constant description 
complexity, is bounded by $O(n^{k-1+\eps})$. Does this bound also hold for
$(T,\pi_1,\pi_2)$-families ? It would be interesting to know if all or
most results in the computational geometry literature relating to arrangements
of sets of constant description complexity, do in fact extend to the more
general setting introduced in this paper. It would also be interesting to
to find proofs of existing bounds using the kind of homological methods
used in this paper. Doing so might remove extraneous assumptions
on general positions in several results  and possibly even lead to 
tighter bounds.

\medskip
\textsl{Acknowledgments}
The author thanks an anonymous referee for several helpful remarks that
helped to substantially improve the paper.

\bibliographystyle{amsplain}
\bibliography{master}

\affiliationone{
   Saugata Basu\\
   Department of Mathematics \\
   Purdue University\\
   West Lafayette, IN 47906 \\
   USA
   \email{sbasu@math.purdue.edu}
}
\end{document}